\documentclass{amsart}
\usepackage[utf8]{inputenc}
\usepackage{amsmath,amssymb}
\usepackage{cite}
\usepackage[final]{showkeys}
\usepackage{hyperref}

\newtheorem{same}{This should never appear}[section]
\newtheorem{defin}[same]{Definition}
\newtheorem{definition}[same]{Definition}

\newtheorem{remark}[same]{Remark}
\newtheorem{theorem}[same]{Theorem}
\newtheorem{example}[same]{Example}

\newtheorem{fact}[same]{Fact}

\newtheorem{cor}[same]{Corollary}

\newtheorem{prop}[same]{Proposition}

\newcommand{\cC}{\mathcal{C}}
\newcommand{\cF}{\mathcal{F}}
\newcommand{\bK}{\mathbb{K}}
\newcommand{\bL}{\mathbb{L}}
\newcommand{\cL}{\mathcal{L}}

\newcommand{\bQ}{\mathbb{Q}}

\newcommand{\ba}{\bold{a}}
\newcommand{\bb}{\bold{b}}
\newcommand{\bc}{\bold{c}}
\newcommand{\bw}{\bold{w}}
\newcommand{\bx}{\bold{x}}
\newcommand{\by}{\bold{y}}
\newcommand{\bz}{\bold{z}}

\newcommand{\cof}{\text{cof}}
\DeclareMathOperator{\Mod}{Mod}
\newcommand{\cf}{\text{cf}}
\newcommand{\rest}{\upharpoonright}
\newcommand{\LS}{\text{LS}}

\newcommand{\ED}{\text{ED}}

\newcommand{\footnotei}[1]{}
\newcommand{\seq}[1]{\langle #1 \rangle}

\title{Cofinality quantifiers in AECs and beyond}
\author{Will Boney}
\date{\today}

\begin{document}

\begin{abstract}
The cofinality quantifiers were introduced by Shelah as an example of a compact logic stronger than first-order logic.  We show that the classes of models axiomatized by these quantifiers can be turned into an Abstract Elementary Class by restricting to \emph{positive} and \emph{deliberate} uses.  Rather than using an ad hoc proof, we give a general framework of abstract Skolemization that can prove a wide range of examples are Abstract Elementary Classes.
\end{abstract}

\maketitle

\section{Introduction}

Abstract Elementary Classes (AECs, introduced by Shelah \cite{sh88}) are the primary framework to do classification theory beyond first-order logic.  They are defined as a collection $(\bK, \prec_\bK)$ of structures $\bK$ and a strong substructure relation $\prec_\bK$ satisfying a certain set of axioms; see \cite{baldwinbook} for an introduction to these axioms and the basic properties of AECs.  These axioms are designed to be broad enough to contain classes axiomatized by $\bL_{\lambda, \omega}$ and some extensions, but provide enough structure to do classification theory.  Beyond AECs, $\mu$-AECs (introduced by Boney, Grossberg, Lieberman, Rosicky, and Vasey \cite{bglrv-muaecs}) can capture classes axiomatized in $\bL_{\lambda, \mu}$, but at the loss of the development of their classification theory.

The primary extension of $\bL_{\lambda, \omega}$ that form AECs are extensions by cardinality quantifiers $Q_\alpha$ (although Baldwin, Ekloff, and Trlifaj \cite{bet-nperp} provide an extension in a different direction).  Here, the cardinality quantifier is interpreted so $Q_\alpha x \phi(x, \by)$ is true (in some structure) iff there are at least $\aleph_\alpha$-many $x$ that make $\phi(x, \by)$ true.  Then classes axiomatized in $\bL_{\lambda, \omega}(Q_\alpha)$ form an AEC, although the strong substructure relation must be strengthened.\footnote{The reason these classes form an AEC and that the substructure relation must be changed is essentially because both $Q_\alpha$ and $\neg Q_\alpha$ can be expressed in an existential fragment of some $\bL_{\lambda, \mu}$; see Section \ref{abskol-sec} for more.}  Most logics extending $\bL_{\lambda, \omega}$ that axiomatize AECs work by adding quantifiers that have a similar `feel' to the cardinality quantifiers, for instance the Ramsey or Magidor-Malitz quantifiers \cite{mm-qmm} or the structure quantifiers \cite{bv-structural}.

We show how extension by another type of quantifier--the cofinality quantifiers introduced by Shelah \cite{sh43} (see Section \ref{back-sec})--can be made an AEC.  Cofinality quantifiers, given a set of regular cardinals $\cC$, $Q^\cof_\cC$ are a binary quantifier, where 
$$Q^\cof_\cC x,y \, \phi(x,y, \bz)$$
means that $\phi(x,y, \bz)$ is a linear order whose cofinality is in $C$.  Among the many properties of cofinality quantifiers, perhaps the most surprising is that $\bL(Q^\cof_\cC)$ is fully compact (see \cite{sh43} or Fact \ref{cof-comp-fact})!  This is a rarity among extensions of first-order logic: Lindstr\"{o}m's Theorem \cite{l-lindstrom} says that any extension of first-order logic must sacrifice either compactness or the downward L\"{o}wenheim-Skolem property, but in practice most extensions of first-order logic become incompact.  This makes cofinality quantifiers particularly intriguing to capture by an AEC since even small fragments of compactness can greatly advance the classification theory (see Boney and Vasey \cite{bv-survey} for a survey).

In order to make cofinality quantifiers an AEC, we must make certain changes the class.  We give more details in Section \ref{aec-sec}, but the essential issue is that the cofinality of a linear order is not preserved under increasing unions.  This necessitates two changes:
\begin{itemize}
    \item Like cardinality quantifiers, the strong substructure relation must preserve the cofinality of linear orders with a positive instance of the cofinality quantifier.  This manifests by not allowing end extensions of such linear orders. 
    \item We have the additional issue that the cofinality can \emph{decrease} following an increasing union.  This requires that we restrict to what we call \emph{positive, deliberate} uses of the cofinality quantifier.  Definition \ref{pos-del-def} makes this precise.
\end{itemize}

Given a positive $\bL(Q^\cof_\cC)$-theory $T$, we form an AEC $\bK_T^+$ through the deliberate use of thises quantifiers (Definition \ref{pos-del-def}), and we briefly explore the properties of this AEC.  An unfortunate consequence of the changes to strong substructure is that many of the nice properties of elementary classes that follow from compactness (amalgamation, etc.) do not hold in these AECs, although these AECs do have some nice properties.  We discuss how some of these issues have their roots in the models produced by the compactness theorem for cofinality quantifiers.  Still, there are some general results that hold for any classes of models that can be made into an AEC with some strong substructure relation (existence of EM models, undefinability of well-order, etc.), and these apply to our classes. Moving beyond AECs, classes axiomatized by cofinality quantifiers naturally form a $\mu$-AEC without the changes above; here, $\mu$ is the successor of the supremum of $C$ for $Q^\cof_\cC$.  This follows from the fact that the cofinality quantifier $Q^\cof_\cC$ is expressible in $\bL_{\infty, \mu}$.

Rather than proving that $\bK_T^+$ forms an AEC through an ad hoc method, we present a general framework of finitary abstract Skolemizations (Definition \ref{abskol-def}).  This captures the essence of Shelah's Presentation Theorem \cite{sh88}, but with a tighter connection to the syntax used to define the AEC.  This is a rather broad method and is able to encompass most known quantifiers that define AECs (see Example \ref{abskol-ex}).

\section{Cofinality quantifiers and background}\label{back-sec}

Background on abstract logics and quantifiers is given in Barwise \cite{b-handbook-math-logic}, but is not really necessary here.  The reader unfamiliar with these ideas can always replace an abstract logic $\cL$ by, depending on the circumstance, one of: finitary first-order logic $\bL = \bL_{\omega, \omega}$; infinitary logic $\bL_{\lambda, \omega}$; or a mild extension of $\bL_{\lambda \omega}$ by cardinality quantifiers.  For completeness, the logic $\bL_{\lambda, \mu}$ (for regular $\lambda \geq \mu$) extends first-order logic by closing formula formation under $<\lambda$-sized disjunctions and conjunctions, and existential and universal quantification of $<\mu$-sized sequences of free variables, with the obvious semantics.

Cofinality quantifiers were introduced by Shelah \cite{sh43} to answer questions of Keisler and Friedman on compact logics stronger than first-order.  We gave an informal description in the Introduction, but give a formal definition here.

\begin{defin}\label{cof-quant-def}
Fix a logic $\cL$, a class of regular\footnote{The cofinality of a linear order is always a regular cardinal, so adding singular cardinals makes no difference} cardinals $\cC$, and a language $\tau$.
\begin{enumerate}
    \item The logic $\cL\left(Q^\cof_\cC\right)$ is an extension of $\cL$ where we add a formulation rule where, if $\phi(\bx, \by, \bz)$ is a formula of $\cL\left(Q^\cof_\cC\right)(\tau)$ (with $\ell(\bx) = \ell(\by)$ finite), then so is 
    $$Q^\cof_\cC \bx,\by \phi(\bx, \by, \bz)$$
    with $\bz$ the remaining free variables.  The semantics of this formula are given by, if $M$ is a $\tau$-structure and $\bc \in M$, then
    $$M\vDash Q^\cof_\cC \bx, \by \phi(\bx, \by, \bc)$$
    iff the relation $\phi(\bx, \by, \bc)$ is a linear order without last element of the set $I:=\{\ba \in M :\text{ there is }\bb \in M, M\vDash \phi(\ba, \bb, \bc) \}$ and that the cofinality of this linear order is in $C$.
    \item A \emph{fragment} $\cF$ of $\cL\left(Q^\cof_\cC\right)(\tau)$ is a collection $\cF \subset \cL\left(Q^\cof_\cC\right)(\tau)$ of formulas that is closed under subformulas.
\end{enumerate}
When we have a singleton $\cC = \{\kappa\}$, we write $Q^\cof_\kappa$ in place of $Q^\cof_{\{\kappa\}}$; there is no risk of confusion because we never place finite cardinals in $\cC$.
\end{defin}

Note that the assertion that `$\phi(\bx, \by, \bz)$ is a linear order without last element' is expressible by a single first-order sentence, and it is the assertion about the cofinality that moves beyond $\bL$. Also, due to the requirement that $\phi(\bx, \by, \bc)$ forms a linear order, we have several equivalent ways to define the set underlying set $I$:
\begin{eqnarray*}
\{\ba \in M :\text{ there is }\bb \in M, M\vDash \phi(\ba, \bb, \bc) \} &=& \{\bb \in M :\text{ there is }\ba \in M, M\vDash \phi(\ba, \bb, \bc) \}\\ &=& \{\ba \in M : M\vDash \phi(\ba, \ba, \bc) \}
\end{eqnarray*}
The last is compactly denoted $\phi(M, M, \bc)$ and is how we will most often refer to the underlying set.

The most common of the cofinality quantifiers used is $Q^\cof_\omega$.  Perhaps the most useful fact about cofinality quantifiers is that first-order logic augmented by a single cofinality quantifier is compact; recall a logic $\cL$ is \emph{compact} iff given any theory $T \subset \cL(\tau)$, $T$ has a model iff every finite subset has a model.

\begin{fact}[{\cite{sh43}, \cite[Corollary 4.4]{cz-cof-comp}}]\label{cof-comp-fact}
For every class $\cC$ of regular cardinals, $\bL\left(Q^\cof_\cC\right)$ is compact.
\end{fact}

\begin{remark}
\begin{enumerate}
    \item In keeping with Lindstr\"{o}m's Theorem, $\bL\left(Q^\cof_\cC\right)$ fails the countable downward L\"{o}wenheim-Skolem property.  For instance, the $\bL\left(Q^\cof_\omega\right)(\{<\})$-sentence
    $$``x < y\text{ is a linear order with no last element''} \wedge \neg Q^\cof_\omega x,y (x < y)$$
    has no countable model.
    \item Casanovas and Ziegler \cite{cz-cof-comp} have recently provided an excellent and self-contained exposition of Fact \ref{cof-comp-fact}. A reader suprised to learn about the full compactness of $\bL(Q^{\cof}_C)$ is not alone; the exposition by Casanovas and Ziegler \cite{cz-cof-comp} was apparently inspired by a referee of Casanovas and Shelah \cite{cs1116} that did not believe this was a ZFC result.
\end{enumerate}
\end{remark}

For future reference, it is helpful to understand the basic structure of the proof of Fact \ref{cof-comp-fact}. At the start two cardinals are fixed: $\kappa \in \cC$ and $\lambda \in \cC$ (if $\cC$ is empty or all regular cardinals, then $Q^\cof_\cC$ is not interesting).  Then $T$ is expanded to definably link all definable linear orders that are positively quantified by $Q^\cof_\cC$ in one group and all definable linear orders that are negatively quantified by $Q^\cof_\cC$ in another group.  Then we find a model of the first-order part of $T$ in which all definable linear orders have cofinality $\max\{\kappa, \lambda\}$.  This is iteratively extended by end-extending all of one group of the linear orders while fixing the other group using the Extended Omitting Types Theorem\cite[Theorem 2.2.19]{changkeisler}.  We continue this iteration for $\min\{\kappa, \lambda\}$-many steps to achieve the desired cofinalities.

The key take away from this proof is that \emph{always} produces models were the definable linear orders have one of exactly two cofinalities.\footnotei{WB: Should I refer to Vaananen's proof?\\ 
\url{http://www.crm.cat/en/Activities/Curs_2016-2017/Documents/Tutorial_CRM_slides.pdf} \\

starting on slide 113, there is a proof of compactness}

Finally, we give an explicit example showing that  $Q^\cof$ is stronger than $\bL_{\infty, \omega}$.

\begin{example}\label{stronger-ex}
Recall that two structures are back and forth equivalent to each other if and only if they are $\bL_{\infty, \omega}$ equivalent.  It is routine to show that $(\bQ, <)$ is aback and forth equivalent to $(\bQ \times \omega_1, <)$.  However, $(\bQ,<)$ satisfies $Q^\cof_\omega x,y (x < y)$ while $(\bQ \times \omega_1, <)$ does not.
\end{example}

We also provide the the basics of AECs (and $\mu$-AECs); \cite{baldwinbook, ramibook} provide further background.

\begin{defin}\label{muaec-def}
Fix an infinite cardinal $\mu$.  A \emph{$\mu$-Abstract Elementary Class} (or $\mu$-AEC for short) is a pair $\left(\bK, \prec_\bK\right)$ where $\bK$ is a collection of structures in a fixed $<\mu$-ary language $\tau_\bK$ satisfying the following axioms
\begin{enumerate}
    \item $\prec_\bK$ is a partial order on $\bK$ that is stronger than $\subset_{\tau_\bK}$.
    \item $\bK$ and $\prec_\bK$ are closed under isomorphisms.
    \item (Coherence) If $M_0, M_1, M_2 \in \bK$ such that $M_0 \prec_\bK M_2$, $M_1 \prec_\bK M_2$, and $M_0 \subset_{\tau_\bK} M_1$, then $M_0 \prec_\bK M_1$
    \item (Closure under $\mu$-directed limits) Given a $\mu$-directed system $\{M_i \in \bK : i \in I\}$, we have that the colimit of this system $\displaystyle\bigcup_{i \in I}M_i$ computed in the category of $\tau$-structures is also the colimit in $\bK$.
    \item (Lowenheim-Sk\"{o}lem-Tarksi number) There\footnote{Formally, once there is a cardinal satisfying this property, all cardinals above it do as well, so we set $\LS(\bK)$ to be the minimal such cardinal.} is a cardinal $\LS(\bK)$ such that, for all $M \in \bK$ and $A \subset M$, there is $M_0 \prec_\bK M$ such that $A \subset M_0$ and $\|M_0\| = |A|^{<\mu}+\LS(\bK)$.
\end{enumerate}
When clear, we often use $\bK$ to refer to the pair $(\bK, \prec_\bK)$.

By far the most common (and important) case is $\mu = \omega$, where we omit $\mu$ and just call it an Abstract Elementary Class (or AEC).
\end{defin}

AECs were introduced by Shelah \cite{sh88}, and generalized to $\mu$-AECs in \cite{bglrv-muaecs}.  They are the most common framework to develop classification theory for nonelementary classes.

\section{$\bL(Q^\cof_\cC)$ as an Abstract Elementary Class}\label{aec-sec}

Fix a set of regular cardinals $\cC$ and a theory $T$ in some fragment $\cF$ of $\bL\left(Q^\cof_\cC\right)(\tau)$.  Recall that a fragment $\cF$ is a subset of $\bL\left(Q^\cof_\cC\right)(\tau)$ that is closed under subformulas.  To build a notion of strong substructure that makes $\Mod T$ an Abstract Elementary Class, we will develop the notion of \emph{positive, deliberate} uses of the cofinality quantifier (Definition \ref{pd-def}).  

The main problem with making $\Mod T$ an AEC is the smoothness under unions of chains: If $\seq{I_\alpha : \alpha < \lambda}$ is a sequence of linear orders such that $I_\alpha$ is end-extended by $I_{\alpha+1}$, then $\displaystyle \bigcup_{\alpha < \lambda} I_\alpha$ has cofinality $\cf \, \lambda$ regardless of the cofinalities of the $I_\alpha$.  This necessitates two changes:
\begin{itemize}
    \item If $M\vDash Q^\cof_\cC x,y \phi(x,y)$, then we should not allow end extensions of this definable linear order in strong extensions; note this is similar to the condition on strong extensions when using the cardinality quantifiers.
    \item If $M \vDash \neg Q^\cof_\cC x,y \phi(x,y)$, then we similarly worry about end extensions.  However, disallowing any end extensions of {\it any} definable linear order would be too restrictive, so we will only allow {\bf positive} instances of the cofinality quantifier.  This doesn't solve the problem completely because definable linear orders that are not put under the $Q^\cof_\cC$ quantifier will `accidentally' end up with a cofinality in $\cC$ after the appropriate unions.  So, via a Morleyization, we avoid this accidental occurence by {\bf deliberately} tagging formulas that we wish to be affected by the cofinality restriction.
\end{itemize}

We detail the construction of positive, definite uses of the cofinality quantifier that will form the strong substructure of an AEC (we deal with positive, deliberate uses of quantifiers in more generality in Section \ref{abskol-sec}).  We work in some degree of generality, allowing for an arbitrary logic $\cL$ to be expanded by cofinality quantifiers, but this will most often be first-order logic $\bL$ with possible extension by infinitary conjunction or cardinality quantifiers.  Note that this expansion is similar to the formation of weak models in \cite{k-existuncount}, but with only one direction of implication.

\begin{defin} \label{pos-del-def}
Fix a language $\tau$ and a logic $\cL$.
\begin{enumerate}
    \item Define $\tau^{\cL}_{*}$ to be 
    $$\tau \cup \{R_\phi(\bz) :\phi(\bx, \by, \bz) \in \cL(\tau)\}$$
    where each $R_\phi$ is new.
    \item Fix a base theory in $\bL\left(Q^\cof_\cC\right)(\tau^\cL_*)$ 
    $$T^{\cof}_{\tau, \cL} := \left\{\forall \bz \left(R_\phi(\bz) \to Q^\cof_\cC \bx, \by \phi(\bx, \by, \bz)\right) : \phi(\bx, \by, \bz) \in \cL(\tau)\right\}$$
    \item Let $T \subset \cL(Q^\cof_\cC)(\tau)$ be a theory where $Q^\cof_\cC$ only appears positively\footnote{The idea of quantifiers appearing positively and the inductive construction of formulas assumed in this definition do not apply to abstract logics generally, but are clearly defined for the logics we will apply this definition to.}.  Define two theories $T^* \subset \cL(\tau^{\cL}_*)$ and $T^+ \subset \cL(Q^{\cof})(\tau^\cL_*)$ by
    \begin{eqnarray*}
    T^* &\text{ is}& \text{ the result of replacing each use of ``$Q^\cof \bx, \by\phi(\bx, \by, \bz)$'' with }\\
    & &\text{``$R_\phi(\bz)$'' in the inductive construction of each $\psi \in T$}\\
    T^+ &:=& T^* \cup T^{\cof}_{\tau,\cL}
\end{eqnarray*}
    \item Given two $\tau^\cL_*$-structures $M \subset N$, define the relation 
    $$M\prec_{\cL(+)} N$$
    iff $M \prec_\cL N$ and, for all $\ba \in M$, if $M \vDash R_\phi(\ba)$, then $\phi(M,M, \ba)$ is cofinal in $\phi(N, N, \ba)$.

\end{enumerate}
\end{defin}

The reason to jump through all of these hoops is the following result.

\begin{theorem}\label{cof-aec-thm}
Fix a set of regular cardinal $C$ and set $\cL = \bL$, finitary first-order logic.  Let $T \subset \cL\left(Q^\cof_\cC\right)(\tau)$ be a theory where all instances of cofinality quantifiers appear positively.  Then
$$\bK_{T}^+ := \left(\Mod T^+, \prec_{\cL(+)}\right)$$
is an Abstract Elementary Class with $\LS(\bK_{T}^+) = |\tau|+(\sup C)^+$.\\

Additionally, the same holds true if $\cL$ is replaced by other logics that axiomatize AECs with the appropriate modification to the strong substructure relation and the L\"{o}wenheim-Skolem number.
\end{theorem}

{\bf Proof:} This is a corollary of the more general result Theorem \ref{nice-quant-thm} \hfill \dag\\

Now that we have an AEC $\bK^+_T$ axiomatized in the fully compact logic $\bL\left(Q^\cof_\cC\right)$, we might hope that several nice consequences of compactness (amalgamation, tameness, etc.) follow directly.  However, this is not the case.  The reason has to do with a disconnect between $\bL(Q^\cof_\cC)$-elementary diagrams and the strong substructure $\prec_{\bK^+_T}$.

Recall that the $\cL$-elementary diagram $\ED_\cL(M)$ of a $\tau$-structure $M$ is the collection of all $\cL\left(\tau \cup \{c_m : m \in M\}\right)$-sententences that are true in $M$ when we interpret $c_m^M = m$.  For first-order logic or any fragment of infinitary $\bL_{\lambda,\kappa}$, we have an equivalence between `there exists an $\cL$-elementary embedding $M \to N$' and `$N \vDash \ED_\cL(M)$.'

However, this does not hold for $\bK^+_T$: modeling $\ED_{\bL\left(Q^\cof_\cC\right)}(M)$ is {\bf not} sufficient to guarantee a $\bK^+_T$-embedding since it does not guarantee that the $R_\phi$-tagged linear orders of $M$ are cofinal in $N$ (Condition \ref{pos-del-def}.(4)).  In the desired uses of compactness, it is routine to find a model of some $\cL$-elementary diagram and turn that into an extension of the desired models; this is not possible.  However, we have some results.

Use $\bL^+(Q^\cof_\cC)$ to denote the $\bL(Q^\cof_\cC)$-formulas where $Q^\cof_\cC$ only appears positively.

\begin{prop}
Let $T$ be an $\bL^+(Q^\cof_\cC)$-theory and define $\bK_T^+$ as in Theorem \ref{cof-aec-thm}.

\begin{enumerate}
	\item $\bK_T^+$ has arbitrarily large models.
	\item Suppose $M \in \bK$ has size $\kappa$ and all $R_\phi$-tagged linear orders have cofinality $\kappa$; that is, if $M \vDash R_\phi(\ba)$, then $\phi(M, M, \ba)$ has cofinality $\kappa$.  Then $M$ has a proper $\prec_{\bK_T^+}$-extension in $\bK_T^+$.
\end{enumerate}
\end{prop}

{\bf Proof:} The first follows easily from the compactness of $\bL(Q^\cof_\cC)$ since it doesn't mention $\prec_{\bK_T^+}$.  For the second, we use the notation and results of \cite{cz-cof-comp}.  Pick some $\psi$ that is not definably connected to the $R_\phi$-tagged linear orders in $M$.  Then, by \cite[Corollary 3.2]{cz-cof-comp}, we can find an $\bL(Q^\cof_\cC)$-elementary extension $N$ of $M$ that extends $\psi$, but in which every $R_\phi$-tagged order has $M$ cofinal.\hfill\dag\\

\subsection{$\bL\left(Q^\cof_\omega\right)$ as an $\omega_1$-Abstract Elementary Class}

Above, there was much effort put into finding precisely the right condition to form an AEC out of $\Mod T$, and the end result was a rather restrictive solution.  Here, we describe a more uniform and natural approach with the drawback that the resulting class is not an Abstract Elementary Class, but instead a $\mu$-Abstract Elementary Class.

While $Q^\cof_\omega$ is not axiomatizable in $\bL_{\infty,\omega}$ (recall Example \ref{stronger-ex}), it is axiomatizable in $\bL_{\omega_1, \omega_1}$.  More generally, for any set of regular cardinals $C$, $Q^\cof_\cC \bx, \by \phi(\bx, \by, \bz)$ is expressible by the first-order statement that $\phi$ defines a linear order without last element and the $\bL_{(\mu+|C|)^+,\mu^+}$ assertion
$$\bigvee_{\lambda \in C} \exists \seq{\bx_i : i < \lambda} \left( \bigwedge_{i<j<\lambda} \phi(\bx_i, \bx_j, \bz) \wedge \forall \bw \bigvee_{i<\lambda} \phi(\bw, \bx_i, \bz)\right)$$
where $\mu = \sup C$.  The logics $\bL_{\lambda,\kappa}$ come with a well-known notion of elementarity.

\begin{definition}
Let $\cC$  be a class of regular cardinals, $\tau$ be a langauge, and $\cF$ be a fragment of $\bL(Q^\cof_C)$.
\begin{enumerate}
    \item Setting $\mu = \sup C$, let $\cF^*$ be the fragment of $\bL_{(\mu+|C|)^+,\mu^+}(\tau)$ that is formed by replacing each instance of $Q^\cof_C$ by the formulation listed above (including the first-order part) and closing under subformula, etc.
    \item Given $\tau$-structures $M$ and $N$, set $M \prec^*_\cF N$ iff $M \prec_{\cF^*} N$.
\end{enumerate}
\end{definition}

\begin{theorem}
Fix a set $C$ of regular cardinals and set $\mu = (\sup C)^+$. For any theory $T$ in $\bL\left(Q^\cof_C\right)(\tau)$, $\bK^*_T$ is a $\mu$-Abstract Elementary Class with $LS(\bK^*_T) = \left(|\tau| + \mu +|C|\right)^{<\mu}$.  This AEC has arbitrarily large models and satisfies the undefinability of well-ordering.
\end{theorem}

\section{Abstract Skolemizations and a sufficient criteria to be an AEC} \label{abskol-sec}

We collect here two related and helpful results: a handy criteria for a class to be an Abstract Elementary Class (Corollary \ref{handy-cor}) and an application of this to show generally that positive, definite uses of infinitary quantification forms an Abstract Elementary Class (Theorem \ref{nice-quant-thm}).

First, we define an abstract notion of Skolemization, that is, an expansion by functions that turns the class into one axiomatizable by a universal theory in $\bL_{\infty,\omega}$.  The goal of this notion is to capture the way in which various extensions of $\bL_{\lambda, \omega}$ by different quantifiers have been turned into AECs.

\begin{defin} \label{abskol-def}
Fix $(\bK, \prec_\bK)$, where $\bK$ is a class of $\tau$-structures and $\prec_\bK$ is a partial order on $\bK$.  A \emph{(finitary) abstract Skolemization (to a universal theory in $\bL_{\infty,\omega}$) of $(\bK, \prec_\bK)$} is an expansion of the langauge $\tau^* := \tau \cup \{F_i : i \in I\}$ by finitary function symbols and a universal theory $T^* \subset \bL_{\infty, \omega}(\tau^*)$ such that the restriction map
$$\cdot \rest \tau :\left(\Mod T^*, \subset\right) \to \left(\bK, \prec_\bK\right)$$
satisfies the following properties:
\begin{enumerate}
    \item (capturing) The restriction map is a functor that is surjective on objects and arrows\footnotei{WB: Is this latter a defined notion? It's what I mean by the first part of what comes next}.
    \item (lifting) Every map $f:M \to N$ in $\bK$ has a lift\footnote{A lift of a model $M$ or an arrow $f:M \to N$ (or a more complicated diagram) from $\bK$ is a model $M^*$ or arrow $f^*:M^* \to N^*$ from $\left(\Mod T^*, \subset\right)$ such that the restriction functor maps them down to the original: $M^* \rest \tau = M$, $N^* \rest \tau = N$, and $f^* \rest \tau = f$} $f^*:M^* \to N*$ in $\Mod T^*$.  Moreover, given any lift $M^*$ of the model $M$ and any map $f:M\to N$ in $\bK$, there is a lift $f^*:M^* \to N^*$ in $\Mod T^*$ with the prescribed domain.
    \item (coherence/local testability) Given $M_0 \subset M_1 \subset N$, if there are separate lifts $M_0^* \subset N^{*0}$ and $M_1^* \subset N^{*1}$, then there are lifts $M_0^{**} \subset M_1^{**} \subset N^{**}$.
\end{enumerate}
We can also define a \emph{$<\mu$-ary abstract Skolemization (to a universal theory in $\bL_{\infty, \mu}$)} by allowing the function symbols to be $<\mu$-ary and the universal theory $T^*$ to be in $\bL_{\infty, \mu}$.\\

We could also speak of abstract Skolemizations to theories in (fragments of) logics \emph{different than} universal theories, but we don't have use for that here.

We often omit `finitary' and `to a universal theory in $\bL_{\infty, \omega}$'.
\end{defin}

We do not explicitly mention it in the definition, but the restriction functor as above is faithful (injective on arrows).

\begin{prop}
Any restriction functor as above is faithful.
\end{prop}

{\bf Proof:} In both categories, the arrows between structures are determined by their value on the underlying sets.\hfill\dag\\

The following sequence of results connects AECs to abstract Skolemizations (one direction is Shelah's Presentation Theorem).

\begin{theorem}\label{one-dir-thm}
If $(\bK, \prec_\bK)$ has an abstract Skolemization, then $(\bK, \prec_\bK)$ is an AEC with L\"{o}wenheim-Skolem number $|\tau^*|$, where $\tau^*$ is the language in the witnessing expansion.
\end{theorem}

{\bf Proof:} Let $\{F_i : i \in I\}$ and $T^* \subset \bL_{\infty, \omega}\left(\tau \cup \{F_i:i\in I\}\right)$ witness the abstract Skolemization.  Most of the AEC axioms (recall Definition \ref{muaec-def} when $\mu=\omega$) follow immediately.  We comment on the three axioms that tend to cause issues for classes being AECs: coherence, smoothness, L\"{o}wenheim-Skolem.

\underline{Coherence:} This is directly addressed by the `coherence/local testability' property of the expansion.  If we have $M_0 \subset M_1$, $M_0 \prec_\bK M_2$, and $M_1 \prec_\bK M_2$, the surjectivity of the restriction gives lifts $M_0^* \subset M_2^{*0}$ and $M_1^* \subset M_2^{*1}$.  This is precisely the set up to give lifts $M_0^{**} \subset M_1^{**} \subset M_2^{**}$.  By applying the restriction functor, we have $M_0 \prec_\bK M_1$, as desired.

\underline{Smoothness:} This is the key use of the condition (2).  Let $\langle M_i \in \bK : i < \alpha\rangle$ be a continuous, $\prec_\bK$-increasing chain of structures with $\alpha$ limit.  We define a continuous, $\subset$-increasing chain $\langle M_i^* \vDash T^* : i < \alpha\rangle$ such that $M_i^*$ is a lift of $M_i$.  To do this, start by letting $M_0^*$ be any lift of $M_0$.  For successors $i = j+1$, we have a lift $M^*_j$ of $M_j$ and $M_j \prec_{\bK} M_i$, so condition (2) guarantees a lift $M^*_i$ of $M_i$ such that $M^*_j \subset M^*_i$.  For limits, we can take unions since the restriction functor preserves unions.\\
In the end, we have that 
$$\bigcup_{i < \alpha} M_i = \left(\bigcup_{i<\alpha} M_i^*\right) \rest \tau$$
so this union is in $\bK$ and is the least upper bound of the chain.

\underline{L\"{o}wenheim-Skolem:} Let $A \subset M \in \bK$ and let $M^*$ be a lift of $M$.  Then, since the restriction functor doesn't change the universe, $A \subset M^*$.  Set $M^*_0$ to be the closure of $A$ under the $\tau\cup\{F_i : i \in I\}$-functions of $M$.  Since $T^*$ is universal, $M^*_0 \vDash T^*$, so $M^*_0 \rest \tau \prec_\bK M$, contains $A$, and has size $\leq |A|+|\tau \cup \{F_i:i\in I\}|$.

\hfill\dag\\

\begin{theorem}\label{other-dir-thm}
If $(\bK, \prec_\bK)$ is an AEC, then the expansion given in Shelah's Presentation Theorem is an abstract Skolemization.
\end{theorem}

The following proof assumes familiarity with the proof and the idea of Shelah's Presentation Theorem; see Baldwin and Boney \cite[Section 3.1]{bb-hanf} for an exposition.

{\bf Proof:} Shelah's Presentation Theorem presents $(\bK, \prec_\bK)$ by an expansion to $\tau^*=\tau(\bK) \cup \{F^n_\alpha : n < \omega, \alpha < \LS(\bK)\}$ that omit a collection $\Gamma$ of quantifier-free types.  We can express this omission through the following $\bL_{\infty, \omega}$ sentence
$$\bigwedge_{p \in \Gamma} \forall \bx \bigvee_{\phi \in p} \neg \phi(\bx)$$
The statement of Shelah's Presentation Theorem (\cite[Fact 3.1.1]{bb-hanf} is perfect for our purposes) gives everything we need except for the coherence/local testability condition.  But this holds exactly because the starting class $(\bK, \prec_\bK)$ is an AEC and, therefore, satisfies coherence.\hfill\dag\\

\begin{cor}\label{handy-cor}
Given a pair $(\bK, \prec_\bK)$ in a finitary language, we have that $(\bK, \prec_\bK)$ is an AEC iff it has a finitary abstract Skolemization to a universal theory in $\bL_{\infty, \omega}$.
\end{cor}

{\bf Proof:} The two directions are Theorems \ref{one-dir-thm} and \ref{other-dir-thm}.\hfill \dag\\

We can also generalize this result to $\mu$-AECs, which we state without proof (the proof is the same).

\begin{theorem}
Given a pair $(\bK, \prec_\bK)$ in a $<\mu$-ary language, we have that $(\bK, \prec_\bK)$ is a $\mu$-AEC iff it has a $<\mu$-ary abstract Skolemization to a universal theory in $\bL_{\infty, \mu}$.
\end{theorem}

Now we turn to the question of how to find abstract Skolemizations.  Shelah's Presentation Theorem gives one way, but the true motivation came from the reason that some quantifiers (like cardinality or cofinality) allow us to form Abstract Elementary Classes: in each case, the quantifiers are expressible in a fragment of $\bL_{\infty, \infty}$ whose only use of inifinitary quantification was a single existential quantifier at the very beginning.
\begin{example}\label{abskol-ex}
\begin{eqnarray*}
Q_\alpha \bx \phi(\bx, \by) &\iff& \exists \langle \bx_i : i < \aleph_\alpha\rangle \left(\bigwedge_{i < \aleph_\alpha} \phi(\bx_i, \by) \wedge \bigwedge_{i<j<\aleph_\alpha} \bx_i \neq \bx_j \right) \\
\neg Q_{\alpha+1} \bx \phi(\bx, \by) &\iff& \exists \langle \bx_i :i < \aleph_\alpha\rangle \forall \bz\left(\phi(\bz, \by) \to \bigvee_{i<\aleph_\alpha} \bz = \bx_i \right)\\
Q^\cof_\kappa \bx, \by \phi(\bx, \by, \bz) &\iff& \exists \langle \bx_i : i < \kappa \rangle \left( `\phi(\bx, \by, \bz) \text{ defines a linear order with no last element'} \right.\\
& & \left. \wedge \forall \bx' \bigvee_{i < \kappa} \phi(\bx', \bx_i, \bz)\right)\\
Q^{ec}_\alpha \bx, \by \phi(\bx, \by, \bz) &\iff& \exists \langle \bx_i : i < \aleph_\alpha\rangle \left( `\phi(\bx, \by, \bz)\text{ defines an equivalence relation} \right.\\
& &\left. \wedge \forall \bx' \left( \phi(\bx', \bx', \bz) \to \bigvee_{i<\aleph_\alpha} \phi(\bx_i, \bx', \bz)\right)\right)  \footnotei{WB: Are there others to add?}
\end{eqnarray*}
Here $Q^{ec}_\alpha$ says that the formula gives a definable equivalence relation with at least $\aleph_\alpha$-many equivalence classes.  The Ramsey/Magidor-Malitz quantifiers \cite{mm-qmm} could be similarly expressed in this way.  Moreover, the strong substructure relation is exactly elementarity according to the fragment of $\bL_{\infty, \infty}$ containing the right hand sides
\end{example}
This form is exactly what allows for a finitary abstract Skolemization to $\bL_{\infty, \omega}$ (which can then be further Skolemized to a universal theory).  Note that the fact that both the positive and negative instances of cardinality quantifiers\footnote{When $\alpha$ is a successor, this is immediate from what is written.  When $\alpha$ is limit, $\neg Q_\alpha x \phi(x)$ is equivalent to $\displaystyle\bigvee_{\beta<\alpha} \neg Q_\beta x \phi(x)$.}have this nice form is what accounts for not needing to worry about deliberate uses of this quantifier.

We now make this connection precise, beginning with some definitions.

\begin{defin}
Say that a quantifier $Q$ is \emph{$\kappa$-existentially definable\footnotei{WB: But there has to be a better name for this concept} in $\cL_1$ over $\cL_0$} iff, for each language $\tau$, there is a map
$$\phi(\bx, \by) \in \cL_0(\tau) \mapsto \Psi_\phi(\bx_i : i < \kappa, \by) \in \cL_1(\tau)$$ such that the following holds
$$\vDash \forall \by \left[\left(Q\bx\phi(\bx, \by)\right) \leftrightarrow \left( \exists \{\bx_i : i < \kappa\} \Psi_\phi(\bx_i : i<\kappa, \by)\right) \right]$$
\end{defin}

The following definition only makes sense due to a subtle (and often overlooked) feature of $\bL_{\infty, \omega}$: the formation of infinitary conjuncts and disjuncts is only allowed if the resulting formula has only finitely many free variables.

\begin{defin}
The logic $\bL_{(\lambda, \omega)}$ is exactly like $\bL_{\lambda, \omega}$ except {\bf without} the restriction to finitely many free variables in conjunctions and disjunctions.
\end{defin}

To emphasize this distinction, consider the formulas around the  assertion that a given $\omega$-sequence is ill-founded
\begin{eqnarray*}
\phi(x_n : n< \omega)&:=& ``\bigwedge_{n < \omega} x_{n+1}<x_n\text{"}\\
\psi(x_n : n < \omega) &:=& ``\phi(x_n:n<\omega) \to \exists y \bigwedge y < x_n\text{"}\\
\Phi &:=& ``\exists \langle x_n : n<\omega\rangle \phi(x_n : n < \omega)\text{"}
\end{eqnarray*}
\begin{itemize}
    \item We have $\phi(\bx), \psi(\bx), \Phi \not\in \bL_{\omega_1, \omega}$ and well-ordering is undefinable in this logic.  Further, we can arrange $M \prec_{\bL_{\omega_1, \omega}} N$ with an ill-founded sequence in $M$ such that a lower bound is added in $N$.
    \item We have $\phi(\bx), \psi(\bx) \in \bL_{(\omega_1,\omega)}$ but $\Phi \not\in \bL_{(\omega_1, \omega)}$, so well-ordering is still not definable in this logic.  However, if $M \prec_{\bL_{(\omega_1, \omega)}} N$ and $\ba=\{a_n \in M : n<\omega\}$ is an ill-founded sequence, then any lower bound in $N$ must already occur in $M$ (by applying the elementarity to $\psi(\ba)$).
    \item We have that $\phi(\bx), \psi(\bx), \Phi \in \bL_{\omega_1, \omega_1}$ and well-ordering is definable in this logic.
\end{itemize}
This definition captures the quantifiers listed above.
\begin{prop}
All quantifiers in Example \ref{abskol-ex} are $\kappa$-existentially definable in $\bL_{(\lambda+\kappa, \omega)}$ over $\bL_{\lambda, \omega}$ for the appropriate $\kappa$.
\end{prop}

{\bf Proof:} The required maps are given in Example \ref{abskol-ex}.\hfill \dag\\

\begin{theorem}\label{nice-quant-thm}
If $Q$ is $\kappa$-existentially definable in $\bL_{(\mu, \omega)}$ over $\bL_{\lambda, \omega}$, then classes axiomatized by positive, deliberate uses of $Q$ in $\bL_{\lambda, \omega}$ have finitary abstract Skolemizations.\\

Furthermore, the same holds if $\bL_{\lambda, \omega}$ is extended by some collection of quantifiers that are $\kappa$-existentially definable in $\bL_{(\mu, \omega)}$ over $\bL_{\lambda, \omega}$.
\end{theorem}

{\bf Proof:} Let $Q$ be $\kappa$-existentially definable in $\bL_{(\lambda, \omega)}$ in $\bL_{(\mu, \omega)}$ over $\bL_{\lambda, \omega}$ via the map $\phi(\bx, \by) \mapsto \Psi(\bx_i : i<\kappa, \by)$.  Following Definition \ref{pos-del-def}, an axiomatization via positive, deliberate uses of $Q$ in $\bL_{\lambda, \omega}$ consists of 
\begin{itemize}
    \item $T \subset \bL_{\lambda, \omega}\left(Q\right)(\tau)$ with $Q$ only occuring positively;
    \item $\tau^* = \tau \cup \{R_\phi(\by) : \phi(\bx, \by) \in \bL_{\lambda, \omega}(\tau)\}$;
    \item $T^* \subset \bL_{\lambda, \omega}(\tau)$ is the result of inductively replacing instances of `$Q\bx \phi(\bx, \by)$' in $T$ by `$R_\phi(\by)$'; and
    \item $T^+ = T^* \cup \left\{\forall \by\left( R_\phi(\by) \to Q \bx \phi(\bx, \by) \right) : \phi(\bx, \by) \in \bL_{\lambda, \omega} \right\}$
\end{itemize}
Then $\bK := \Mod(T^+)$ is the class we need to we need to provide the Skolemization for.  We describe the Skolemization in two steps.

For the first step, for each $\phi(\bx, \by)$, add functions 
$$ \left\{ F^{Q, \phi}_{i, j}(\by) : i < \kappa, j < \ell(\bx)\right\}$$
and set
$$T^{++} = T^* \cup \left\{ \forall \by\left( R_\phi(\by) \to \Psi_\phi\left(F^{Q, \phi}_{i,j}(\by) : j < \ell(\bx), i<\kappa, \by\right)\right) : \phi(\bx, \by) \in \bL_{\lambda, \omega} \right\}$$
Crucially, $T^{++}$ is an $\bL_{\lambda+\mu,\omega}$-theory.  So this gives a Skolemization of $T^+$ to a (non-universal) theory in $\bL_{\lambda+\mu, \omega}$.  It is a standard result (e.g., see \cite[Theorem 17]{k-lw1w} for the case $\lambda+\mu=\omega_1$) that $\bL_{\infty, \omega}$ theories have finitary Skolemizations to universal theories in $\bL_{\infty, \omega}$; the second step is to do this Skolemization.

Putting these steps together, we have a finitary Skolemization of $\bK$ to a universal theory in $\bL_{\infty, \omega}$; we can define $\prec_\bK$ by setting $M \prec_\bK N$ iff there are lifts $M^*$ and $N^*$ such that $M^* \subset N^*$. \hfill \dag\\

\begin{cor}
All of the quantifiers listed in Example \ref{abskol-ex} form AECs when used positively and delibarately over $\bL_{\infty, \omega}$ and can be mixed together.
\end{cor}

Note that the strong substructure relation $\prec_\bK$ in Defintion \ref{abskol-def} can be recovered from the expansion $T^*$.  Chasing through the definitions, the appropriate strong substructure relation in the cases above is elementarity according to the fragment of $\bL_{\infty, \infty}$ needed to define the quantifiers; this corresponds exactly to the seemingly ad hoc notions given for cardinality and cofinality quantifiers.

\bibliographystyle{amsalpha}
\bibliography{bib}

\newcommand{\etalchar}[1]{$^{#1}$}
\providecommand{\bysame}{\leavevmode\hbox to3em{\hrulefill}\thinspace}
\providecommand{\MR}{\relax\ifhmode\unskip\space\fi MR }
\providecommand{\MRhref}[2]{%
  \href{http://www.ams.org/mathscinet-getitem?mr=#1}{#2}
}
\providecommand{\href}[2]{#2}
\begin{thebibliography}{BGL{\etalchar{+}}16}

\bibitem[Bal09]{baldwinbook}
John Baldwin, \emph{Categoricity}, University Lecture Series, no.~50, American
  Mathematical Society, 2009.

\bibitem[Bar82]{b-handbook-math-logic}
Jon Barwise (ed.), \emph{Handbook of mathematical logic}, 1st ed.,
  North-Holland Publishing Company, 1982.

\bibitem[BB17]{bb-hanf}
Will Boney and John Baldwin, \emph{Hanf numbers and presentation theorems in
  {AEC}s}, Beyond First-Order Model Theory (Jose Iovino, ed.), CRC Press, 2017,
  pp.~327--352.

\bibitem[BET07]{bet-nperp}
John Baldwin, Paul Eklof, and Jan Trlifaj, \emph{${}^\prep n$ as an abstract
  elementary class}, Annals of Pure and Applied Logic \textbf{149} (2007),
  no.~1-3, 25--39.

\bibitem[BGL{\etalchar{+}}16]{bglrv-muaecs}
Will Boney, Rami Grossberg, Michael Lieberman, Ji{\v r}{\'\i} Rosick{\'y}, and
  Sebastien Vasey, \emph{$\mu$-{A}bstract {E}lementary {C}lasses and other
  generalizations}, Journal of Pure and Applied Algebra \textbf{220} (2016),
  no.~9, 3048--3066.

\bibitem[BV17]{bv-survey}
Will Boney and Sebastien Vasey, \emph{A survey of tame {A}bstract {E}lementary
  {C}lasses}, Beyond First-Order Model Theory (Jose Iovino, ed.), CRC Press,
  2017, pp.~353--.

\bibitem[BV19]{bv-structural}
\bysame, \emph{Structural logic and abstract elementary classes with
  intersections}, Bulletin of the Polish Academy of Sciences Mathematics
  \textbf{67} (2019), no.~1, 1--17.

\bibitem[CK12]{changkeisler}
C.~C. Chang and H.~Jerome Keisler, \emph{Model theory}, 3rd ed., Dover
  Publications, 2012.

\bibitem[CS19]{cs1116}
Enrique Casanovas and Saharon Shelah, \emph{Universal theories and compactly
  expandable models}, Journal of Symbolic Logic \textbf{84} (2019), no.~3,
  1215--1223.

\bibitem[CZ20]{cz-cof-comp}
Enrique Casanovas and Martin Ziegler, \emph{An exposition of the compactness of
  $l(q^{cf})$}, Bulletin of Symbolic Logic \textbf{26} (2020), no.~3-4,
  212--218.

\bibitem[Gro1X]{ramibook}
Rami Grossberg, \emph{A course in model theory}, In preparation, 201X.

\bibitem[Kei70]{k-existuncount}
H.~Jerome Keisler, \emph{Logic with the quantifier ``there exist uncountably
  many''}, Annals of Mathematical Logic \textbf{1} (1970), no.~1, 1--93.

\bibitem[Kei71]{k-lw1w}
\bysame, \emph{Model theory for infinitary logic}, Studies in Logic and the
  Foundations of Mathematics, vol.~62, North-Holland Publishing Company, 1971.

\bibitem[Lin69]{l-lindstrom}
Per Lindstr{\"o}m, \emph{On extensions of elementary logic}, Theoria
  \textbf{35} (1969), no.~1, 1--11.

\bibitem[MM77]{mm-qmm}
Menachem Magidor and Jerome Malitz, \emph{Compact extensions of {$L(Q)$}, part
  1a}, Annals of Mathematical Logic \textbf{11} (1977), 217--261.

\bibitem[She75]{sh43}
Saharon Shelah, \emph{Generalized quantifiers and compact logic}, Transactions
  of the American Mathematical Society \textbf{204} (1975), 342--364.

\bibitem[She87]{sh88}
\bysame, \emph{Classification of nonelementary classes, {II}. {A}bstract
  elementary classes}, Classification theory (John Baldwin, ed.), 1987,
  pp.~419--197.

\end{thebibliography}

\end{document}